\documentclass{amsart}

\title{Uniqueness of Morava $K$-theory}
\author{Vigleik Angeltveit}

\usepackage{amsxtra}
\usepackage{amsfonts}
\usepackage[latin1]{inputenc}
\usepackage{graphicx}
\usepackage{amsmath,amssymb,latexsym,amsthm,mathrsfs}
\usepackage[all]{xy}

\newtheorem{theorem}{Theorem}[section]
\newtheorem{thm}[theorem]{Theorem}
\newtheorem{lemma}[theorem]{Lemma}

\newtheorem{prop}[theorem]{Proposition}

\newtheorem{conj}[theorem]{Conjecture}

\newtheorem{mainthm}{Theorem}

  \newcommand{\cC}{\mathcal{C}} \newcommand{\cD}{\mathcal{D}}       \newcommand{\cK}{\mathcal{K}}  \newcommand{\cM}{\mathcal{M}}

\newcommand{\sA}{\mathscr{A}}

     \newcommand{\bF}{\mathbb{F}}          
         \newcommand{\bZ}{\mathbb{Z}}

\newcommand{\sma}{\wedge} 

\newcommand{\hE}[1]{\widehat{E(#1)}}
\newcommand{\BP}[1]{BP \langle #1 \rangle}
\newcommand{\BR}{B\mathscr{A}}
\newcommand{\ttau}{\bar{\tau}}
\newcommand{\xxi}{\bar{\xi}}
\newcommand{\aalpha}{\tilde{\alpha}}
\newcommand{\tq}{\tilde{q}}

\begin{document}

\begin{abstract}
We show that there is an essentially unique $S$-algebra structure on the Morava $K$-theory spectrum $K(n)$, while $K(n)$ has uncountably many $MU$ or $\hE{n}$-algebra structures. Here $\hE{n}$ is the $K(n)$-localized Johnson-Wilson spectrum. To prove this we set up a spectral sequence computing the homotopy groups of the moduli space of $A_\infty$ structures on a spectrum, and use the theory of $S$-algebra $k$-invariants for connective $S$-algebras found in \cite{DuSh} to show that all the uniqueness obstructions are hit by differentials.
\end{abstract}

\address{Department of Mathematics \\ University of Chicago \\ Chicago IL 60637}


\maketitle

\section{Introduction}
We study the moduli space of $S$-algebra structures on the Morava $K$-theory spectrum $K(n)$. Recall that given a prime $p$ and an integer $n \geq 1$, $K(n)$ is the spectrum carrying the Honda formal group of height $n$ over $\bF_p$, and that $K(n)_* \cong \bF_p[v_n,v_n^{-1}]$ with $|v_n|=2p^n-2$. Robinson \cite{Ro88} found that there are uncountably many ways to build an $A_\infty$ structure on $K(n)$, but he did not ask if these $A_\infty$ structures might all be equivalent. The point is that there are two distinct definitions of the moduli space of $S$-algebra structures, and in this paper we use the version where we allow automorphisms of the underlying spectrum. We prove the following:

\begin{mainthm} \label{thm:main1}
There is an essentially unique $S$-algebra structure on $K(n)$, in the sense that the moduli space of $S$-algebra structures on $K(n)$ is connected.
\end{mainthm}

This should be compared to the situation where we study the moduli space of $R$-algebra structures on $K(n)$ for some other commutative $S$-algebra $R$:

\begin{mainthm} \label{thm:main2}
Let $R=MU$ or $R=\hE{n}$. Then there are uncountably many $R$-algebra structures on $K(n)$, in the sense that the moduli space of $R$-algebra structures on $K(n)$ has uncountably many path components.
\end{mainthm}
\noindent
If $BP$ is a commutative $S$-algebra, Theorem \ref{thm:main2} remains true with $R=BP$.

We will use two approaches to study the moduli space of $S$-algebra or $R$-algebra structures on a spectrum $A$. For our first approach, we use the equivalence between $S$-algebras and $A_\infty$ ring spectra, and study how to build an $A_\infty$ structure on $A$ by induction on the $A_m$ structure. This is the approach taken by Robinson \cite{Ro88} and later by the author \cite{AnTHH}. We need to modify this approach slightly to get the right notion of equivalence of $A_\infty$ structures; this amounts to allowing maps $(A,\phi) \to (A,\psi)$ of $A_\infty$ ring spectra where the underlying map $A \to A$ of spectra is not the identity but merely a weak equivalence.

We will define the appropriate moduli space of $S$-algebra structures on $A$, which we denote by $\BR^S(A)$, and set up a spectral sequence $\{E_r^{s,t}\}$ which contains the obstructions to $\BR^S(A)$ being nonempty and, given a basepoint, computes the homotopy groups of this space. The spectral sequence is similar to the one found in \cite{Re97} based on derived functors of derivations.

Using this approach, the uniqueness obstructions for $K(n)$ lie in $E_\infty^{s,s}$ for $s \geq 1$. On the $E_2$ term, everything in positive filtration is concentrated in even total degree, so every class in $E_2^{s,s}$ for $s \geq 1$ is a permanent cycle. But $E_1^{0,1}$ is large, in fact $E_1^{0,1}$ is closely related to the Morava stabilizer group, and there are potential differentials $d_r : E_r^{0,1} \to E_r^{r,r}$ for all $r \geq 1$ killing the uniqueness obstructions.

This should be compared to the situation for the Morava $E$-theory spectrum $E_n$. See \cite{Re97} for a spectral sequence which computes the space of $A_\infty$ structures on $E_n$ and \cite{GoHo} for a spectral sequence which computes the space of $E_\infty$ structures on $E_n$. In both cases the $E_2$ term is trivial in positive filtration, so there is no need to compute any differentials.

The other approach, which works only if $A$ is connective, is to study how to build $A$ as a Postnikov tower in the category of $S$-algebras. For this we use a result of Dugger and Shipley \cite{DuSh} which tells us that the set of ways to build $P_m A$ from $P_{m-1} A$ as an $S$-algebra can be calculated using $THH_S^{m+2}(P_{m-1} A; H\pi_m A)$.

These topological Hochschild cohomology groups can be calculated when $A=k(n)$ is connective Morava $K$-theory, and this lets us identify the uniqueness obstructions for building $k(n)$ as an $S$-algebra. Once again the obstructions are nontrivial, but something interesting happens. Each of the obstructions we found using the first approach also live in the $E_2$ term of the canonical spectral sequence converging to $THH_S^*(P_{m-1} k(n); H\bF_p)$ for some $m$, but in every case the obstruction is killed by a differential. Hence the corresponding $S$-algebra structures on $P_m k(n)$ are equivalent, and this equivalence can be lifted first to $k(n)$ and then to $K(n)$.

We emphasize that both approaches are necessary to prove Theorem \ref{thm:main1}. Using only the first approach is insufficient because we do not know how to calculate the differentials in the spectral sequence converging to $\pi_* \BR^S(K(n))$ directly. Using only the second approach is insufficient because the connective Morava $K$-theory spectrum $k(n)$ does \emph{not} have a unique $S$-algebra structure. While the obstructions we found in the first approach are killed in the spectral sequence converging to $THH_S^*(P_{m-1} k(n);H\bF_p)$ for suitable $m$, there are other uniqueness obstructions here and we do not have a direct way to show that those obstructions become trivial when inverting $v_n$.

\subsection{Organization}
In \S \ref{s:moduli_space} we define the moduli space of $A_\infty$ structures on a spectrum $A$ and construct a spectral sequence converging to the homotopy groups of this moduli space. Because we need to allow maps of $A_\infty$ ring spectra which commute with the operad structure only up to homotopy and higher homotopies, we use a certain multicategory with $r$ colors to define $(r-1)$-fold composites, and as a result the moduli space is (the geometric realization of) an $\infty$-category, regarded as a simplicial set.

In \S \ref{s:BKSSKn} we compute the $E_2$ term of this spectral sequence for $K(n)$, with $\hE{n}$, $MU$ and $S$ as the ground ring in \S \ref{ss:BKSShEn}, \ref{ss:BKSSMU} and \ref{ss:BKSSS} respectively. In the first two cases the spectral sequence collapses at the $E_2$ term, in the last case there are potential differentials. Counting the classes that are left in $E_2^{s,s}$ with $\hE{n}$ or $MU$ as the ground ring then proves Theorem \ref{thm:main2}.

In \S \ref{s:k_invariants} we recall the theory of $k$-invariants for connective $S$-algebras, which live in topological Hochschild cohomology, due to Dugger and Shipley \cite{DuSh} and discuss the relationship with additive $k$-invariants.

In \S \ref{s:kforkn} we compute the relevant topological Hochschild cohomology groups for Postnikov sections of connective Morava $K$-theory, with $\BP{n}_p$, $MU$ and $S$ as the ground ring in \S \ref{ss:k_invBPn}, \ref{ss:k_invMU} and \ref{ss:k_invS} respectively. The calculation with $\BP{n}_p$ as the ground ring requires optimistic assumptions about the commutativity of the multiplication on $\BP{n}_p$, we include it because it is parallel to the situation of $K(n)$ as an $\hE{n}$-algebra and it gives a clearer conceptual picture of what is going on.

In \S \ref{s:proof} we put the pieces together to prove Theorem \ref{thm:main1}.

Finally, in \S \ref{s:2periodic} we discuss the moduli space of $S$-algebra structures on the $2$-periodic version $K_n$ of Morava $K$-theory. In this case we do not have a unique $S$-algebra structure on $K_n$, but we conjecture that there are only finitely many.

\subsection{Acknowledgements}
The author would like to thank Haynes Miller for his encouragement and for reading an earlier draft of this paper, and Jacob Lurie for help with Lemma \ref{lem:Jacob}.

\section{The moduli space of $A_\infty$ structures} \label{s:moduli_space}
Recall that in a good category of spectra, such as \cite{EKMM}, any $A_\infty$ ring spectrum can be replaced with a weakly equivalent $S$-algebra. Moreover, the functor from the multicategory describing $n$-fold composition of $A_\infty$ ring spectra to the multicategory describing $n$-fold composition of $S$-algebras is a weak equivalence on all $Hom$ sets, and this implies that the moduli space of $A_\infty$ structures on $A$ we define below, which only depends on the homotopy type of $A$, is weakly equivalent to the moduli space of $S$-algebra structures on $A$.


Other approaches to studying the moduli space of $A_\infty$ structures on a spectrum $A$, such as the one found in \cite{Re97}, assumes that $A$ comes with a fixed homotopy commutative multiplication. At $p=2$ the Morava $K$-theory spectrum $K(n)$ does not have a homotopy commutative multiplication \cite{Na02}, and in any case we prefer to fix as little data as possible, so instead of following \cite{Re97} we will set up a similar spectral sequence based on the obstruction theory in \cite{Ro88} and \cite{AnTHH}.

We take an $A_\infty$ ring spectrum to mean an algebra over the Stasheff associahedra operad $\cK=\{K_n\}_{n \geq 0}$. For $0 \leq n \leq \infty$ an $A_n$ structure on $X$ is a compatible family of maps
\[ (K_m)_+ \sma X^{(m)} \to X \]
for $m \leq n$, where $X^{(m)}$ denotes the $m$-fold smash product of $X$ with itself. If we work in the category of $R$-modules for some commutative $S$-algebra $R$, all smash products are over $R$.

Using only maps $X \to Y$ of $A_n$ ring spectra which commute strictly with the operad action is too restrictive, so following Boardman and Vogt \cite{BoVo73} we define a map of $A_n$ ring spectra to be a family of maps $(L_m)_+ \sma X^{(m)} \to Y$ for $m \leq n$, where $L_m$ is a certain polyhedron of dimesnions $m-1$. Here $L_m$ can be defined in terms of the $W$-construction on the multicategory (colored operad, colored PRO) with two objects $0$ and $1$ and $Hom(\epsilon_1,\ldots,\epsilon_n;\epsilon)$ a point if $\epsilon_1+\ldots+\epsilon_n \leq \epsilon$ and empty otherwise, or more concretely as a certain space of metric trees with two colors. We think of $(L_m)_+ \sma X^{(m)} \to Y$ as a homotopy between the maps $(K_m)_+ \sma X^{(m)} \to X \to Y$ and $(K_m)_+ \sma X^{(m)} \to (K_m)_+ \sma Y^{(m)} \to Y$.

As observed in \cite[Ch.\ 4]{BoVo73}, while it is possible to ``compose'' the maps we just defined, composition is not associative. Instead, we get an $\infty$-category (quasi-category, restricted Kan complex) of $A_n$ ring spectra encoding the various ways of composing multiple maps, where an $r$-simplex is a ``composite of $r-1$ maps'' defined in terms of a multicategory with $r$ colors. This is not actually a problem for us, because we can take the geometric realization of an $\infty$-category just as easily as we can take the geometric realization of (the nerve of) a category.

If $R$ is a commutative $S$-algebra and $A$ is an $R$-module, let $\BR^R_n(A)$ be the moduli space of $A_n$ structures on $A$ in the category of $R$-modules. To be precise, we let $\BR^R_n(A)$ be the geometric realization of the $\infty$-category $\sA_n^R(A)$ defined as follows. An object ($0$-simplex) in $\sA_n^R(A)$ is a pair $(X, \phi)$ where $X$ is weakly equivalent to $A$ and $\phi=\{\phi_m\}_{0 \leq m \leq n}$ is an $A_n$ structure on $X$ in the category of $R$-modules. For convenience we will assume $X$ is cofibrant as an $R$-module. A morphism ($1$-simplex) $(X,\phi) \to (Y,\psi)$ is a map $X \to Y$ of $A_n$ ring spectra, where the underlying map $X \to Y$ of spectra is a weak equivalence. An $r$-simplex is defined similarly, as in \cite[Definition 4.7]{BoVo73}. A choice of weak equivalence $X \to A$ is not part of the data. Some care is needed to make sure that we end up with a small ($\infty$-)category, which we need to apply geometric realization, we refer the reader to \cite{DwKa84} for one possible solution.

A general argument due to Dwyer and Kan \cite{DwKa84} shows that the moduli space $\BR_n^R(A)$ decomposes as
\[ \BR_n^R(A) \simeq \coprod_{[X]} BAut_{\sA_n^R(A)}(X), \]
where the coproduct runs over one representative from each path component of $\BR_n^R(A)$ and $Aut_{\sA_n^R}(X)$ is the topological monoid of self-equivalences of a cofibrant-fibrant model for $X$.

In particular, an $A_1$ structure consists only of the unit map $R \to A$ and the identity map $A \to A$, and an automorphism of $A$ as an $A_1$-algebra is a unit-preserving weak equivalence $A \to A$ of $R$-modules. Let $Aut_R(A)_1$ denote the space of unit-preserving $R$-module automorphisms of a cofibrant-fibrant model of $A$. Then $\BR_1^R(A) \simeq B Aut_R(A)_1$.

Given a tower of fibrations
\[ \ldots \to X_n \to X_{n-1} \to \ldots \to X_0\]
with inverse limit $X$, recall \cite[Ch IX, \S 4]{BoKa74} that we get a ``fringed'' spectral sequence (called ``the (extended) homotopy spectral sequence'' in loc.\ cit.)
\[ E_1^{s,t}=\pi_{t-s} F_s \Longrightarrow \pi_{t-s} X, \]
where $F_s$ is the fiber of $X_s \to X_{s-1}$. This is not quite a spectral sequence in the usual sense, for the following reasons. First, $X$ might be empty, and the spectral sequence only exists as long as we can lift a given basepoint up the tower. The terms $E_1^{s,s+1}$ on the superdiagonal, contributing to $\pi_1 X$, are in general nonabelian, and the terms $E_1^{s,s}$ on the diagonal, contributing to $\pi_0 X$, are only sets. The fringing refers to the lack of negative dimensional terms to receive differentials.

This spectral sequence has good convergence properties, it converges completely as long as there are no $\lim^1$ terms \cite[Lemma IX.5.4]{BoKa74}.

Also recall \cite{Bo89} that if the tower of fibrations comes from the Tot-tower of a (simple, fibrant) cosimplicial space, the above spectral sequence has (some) negative dimensional terms. In particular $E_1^{s,s-1}$ exists and serves both as the target of differentials from the diagonal and as the place where obstructions to lifting a basepoint up the tower lie.

In our case the $n$'th space in the tower of fibrations will be the space $\BR^R_{n+1}(A)$, and although this tower does not come from a cosimplicial space we will describe sets $E_1^{s,s-1}$ containing the obstructions to lifting a basepoint up the tower. Moreover, the only nonabelian group on the superdiagonal is $E_1^{0,1}$ and while $E_1^{s,s}$ is not a group, it is a torsor over an abelian group that can be described in the same way as $E_1^{s,t}$ for $t-s \geq 1$.

We wish to identify the fiber of $\BR^R_{n+1}(A) \to \BR^R_n(A)$ with the space of extensions of a given $A_n$ structure on $A$ to an $A_{n+1}$ structure. If $\BR^R_n(A)$ was the classifying space of a category we could use Quillen's Theorem B \cite{Qu73}. Instead we use the following version, with notation from \cite{Lu08}:

\begin{lemma} \label{lem:Jacob}
Suppose $F : \cC \to \cD$ is a map of $\infty$-categories with the property that for every $f : d \to d'$ in $\cD$ the maps
\[ \cC \times_{\cD} \cD_{d/} \overset{\simeq}{\leftarrow} \cC \times_{\cD} \cD_{f/} \overset{\simeq}{\to} \cC \times_{\cD} \cD_{d'/} \]
are weak equivalences. Then the homotopy fiber of $\cC \to \cD$ is weakly equivalent to $\cC \times_{\cD} \cD_{d/}$.
\end{lemma}

\begin{proof}[Sketch proof]
The homotopy fiber of $\cC \to \cD$ is the fiber of
\[ p : \cC \times_{\cD} Fun(\Delta^1,\cD) \to \cD. \]
The hypothesis imply that the inverse image of any $0$-simplex or $1$-simplex in $\cD$ is weakly equivalent to $\cC \times_{\cD} \cD_{d/}$, and the case for a general simplex in $\cD$ follows.
\end{proof}

Let $\bar{A}$ denote the cofiber of the unit map $R \to A$ (assuming $A$ is cofibrant) and let $\bigvee^n \! A$ denote the ``fat wedge''
\[\bigvee\! {}^n \! A = \bigvee_{1 \leq i \leq n} A^{(i-1)} \sma R \sma A^{(n-i)}.\]
Then the canonical map $\bigvee^n \! A \rightarrow A^{(n)}$ is a cofibration, with cofiber $\bar{A}^{(n)}$.

Now consider the forgetful functor $F : \sA_{n+1}^R(A) \to \sA_n^R(A)$. Given $(X,\phi) \in \sA_n^R(A)$ and $(Y,\psi) \in \sA_n^R(A)_{(X,\phi)/}$, the fiber over 
\[ (Y,\psi) \in \sA_{n+1}^R(A) \times_{\sA_n^R(A)} \sA_n^R(A)_{(X,\phi)/} \]
is the space of extensions of the $A_n$ structure $\psi$ on $Y$ to an $A_{n+1}$ structure.

An $A_{n+1}$ structure on $Y$ extending $\psi$ is a map
\[ m_{n+1} : (K_{n+1})_+ \sma Y^{(n+1)} \rightarrow Y \]
satisfying two conditions. First, $m_{n+1}$ is determined by $\psi$ on $(\partial K_{n+1})_+ \sma Y^{(n+1)}$, and second, $m_{n+1}$ is determined by the unitality condition on $(K_{n+1})_+ \sma \bigvee^{n+1} Y$.

The cofiber of the map
\[ (\partial K_{n+1})_+ \sma Y^{(n+1)} \!\!\!\!\!\!\!\! \coprod_{(\partial K_{n+1})_+ \sma \bigvee^{n+1} \! Y} \!\!\!\!\!\!\!\! (K_{n+1})_+ \sma \bigvee \! {}^{n+1} Y \rightarrow (K_{n+1})_+ \sma Y^{(n+1)} \]
is $\Sigma^{n-1} \bar{Y}^{n+1}$, and hence the space of extensions of $\psi$ to an $A_{n+1}$ structure is weakly equivalent to $Hom(\Sigma^{n-1} \bar{Y}^{(n+1)}, Y)$, which is weakly equivalent to
\[ Hom(\Sigma^{n-1} \bar{A}^{(n+1)}, A).\]

Similarly, given $f : (X,\phi) \to (Y,\psi)$ in $\sA_n^R(A)$ and an element $(Z,\xi) \in \sA_n^R(A)_{f/}$, the fiber over $(Z,\xi)$ is the space of extensions of the $A_n$ strucutre $\xi$ on $Z$, and the maps in the Lemma \ref{lem:Jacob} are clearly weak equivalences. Hence we can conclude that the fiber of $F : \BR_{n+1}^R(A) \to \BR_n^R(A)$ is the space of extensions of a given $A_n$ structure to an $A_{n+1}$ structure, as we wanted.

\begin{thm}
There is a spectral sequence $\{E_r^{s,t}\}$ with $E_1^{s,t}$ defined for $s \geq 0$ and $t-s \geq -1$ converging to $\pi_{t-s} \BR^R(A)$ with the obstructions to $\BR^R(A)$ being nonempty on the subdiagonal $t-s=-1$. We have $E_1^{0,-1}=\varnothing$, $E_1^{0,0}=0$, $E_1^{0,1} \cong \pi_0 Aut_R(A)_1$ and
\[ E_1^{s,t} \cong [\Sigma^{t-1} \bar{A}^{(s+1)},A] \]
otherwise. Here $E_1^{s,t}$ is a group for $t-s \geq 1$, a torsor over the corresponding group for $t-s=0$, and a set for $t-s=-1$. 
\end{thm}

\begin{proof}
From the obstruction theory developed in \cite{AnTHH} we conclude that we get a tower of fibrations
\begin{equation*}
\BR^R(A) \simeq \BR_\infty^R(A) \to \ldots \to \BR^R_2(A) \to \BR^R_1(A),
\end{equation*}
and the spectral sequence is the one associated with this tower.

The above discussion identifies $E_1^{s,t}$ for $t-s \geq 0$. The obstruction theory in \cite{AnTHH} also identifies the obstruction to extending an $A_n$ structure to an $A_{n+1}$ structure with an element in $E_1^{n,n-1}$.
\end{proof}

We would like to compare this to topological Hochschild cohomology, in particular to the $E_2$ term of the topological Hochschild cohomology spectral sequence, because that is something we can compute. Let $\{\tilde{E}_r^{p,q}\}$ be the spectral sequence with $E_1$-term $\tilde{E}_1^{p,q}=\pi_q F_S(\bar{A}^{(p)},A)$ converging to $\pi_{q-p}THH_R(A)$ if $A$ is an $R$-algebra so that topological Hochschild cohomology is defined.

\begin{thm}
Suppose $\BR^R(A)$ is nonempty, and choose an $R$-algebra structure on $A$. Then
\[ E_2^{s,t} \cong \tilde{E}_2^{s+1,t-1} \]
for $s \geq 2$ and $t-s\geq -1$. This isomorphism of $E_2$ terms is an isomorphism of abelian groups for $t-s \geq 1$, of torsors for $t-s=0$ and of sets for $t-s=-1$.
\end{thm}

\begin{proof}
The $E_1$-terms are isomorphic for $s \geq 1$ and $t-s \geq -1$, and the argument for why the $d_1$ differential on $E_1^{*,*}$ is isomorphic to the Hochschild differential is contained in \cite{Ro88} or \cite{AnTHH}.
\end{proof}

\section{The spectral sequence for Morava $K$-theory} \label{s:BKSSKn}
In this section we prove Theorem \ref{thm:main2} by explicitly calculating the $E_\infty$ term of the spectral sequence converging to $\pi_* \BR^R(K(n))$ for $R=\hE{n}$ and $R=MU$. We also calculate the $E_2$ term for $R=S$.

\subsection{Ground ring $R=\hE{n}$} \label{ss:BKSShEn}
Let $R=\hE{n}$ be the $K(n)$-localization of the Johnson-Wilson spectrum, with homotopy groups
\[ \hE{n}_* = \bZ_{(p)}[v_1,\ldots,v_{n-1}, v_n^{\pm 1}]^\wedge_I.\]
Here $I=(p,v_1,\ldots,v_{n-1})$ and $(-)^\wedge_I$ denotes $I$-completion. Then $\hE{n}$ can be given the structure of a commutative $S$-algebra \cite{RoWh02}, and $K(n) \simeq \hE{n}/I$. As in \cite{AnTHH} we find that the spectral sequence converging to $\pi_*THH_{\hE{n}}(K(n))=THH^{-*}_{\hE{n}}(K(n))$ collapses at the $E_2$ term (there are interesting extensions) because everything is concentrated in even total degree, with
\[ \tilde{E}_2=\tilde{E}_\infty = K(n)_*[q_0,\ldots,q_{n-1}]. \]
Here $q_i$ is in filtration $1$ and total homological degree $-2p^i$. There can obviously be no $\lim^1$ terms, so the spectral sequence converges completely.

This gives us the positive filtration part of the spectral sequence converging to $\pi_* \BR^{\hE{n}}(K(n))$. In particular, there are no obstructions to the existence of an $A_\infty$ structure on $K(n)$, and the part contributing to $\pi_0 \BR^{\hE{n}}(K(n))$ is the homological degree $-2$ part of $K(n)_*[q_0,\ldots,q_{n-1}]$ of degree at least $2$ in the $q_i$'s. We also know that
\[ \pi_* F_{\hE{n}}(K(n),K(n)) \cong \Lambda_{K(n)_*}(Q_0,\ldots,Q_{n-1}), \]
where $Q_i$ is the Bockstein corresponding to $v_i$, and the degree $0$ part of this is $\bF_p$. Only one of these $p$ maps commutes with the unit $\hE{n} \to K(n)$, so we find that $E_1^{0,1}=E_2^{0,1}=0$. Hence there are no possible differentials originating from $E_2^{0,1}$. Everything in positive filtration is concentrated in even total degree, so the spectral sequence collapses at the $E_2$ term with infinitely many classes on the diagonal. This proves Theorem \ref{thm:main2} for $R=\hE{n}$.

\subsection{Ground ring $R=MU$} \label{ss:BKSSMU}
A similar argument shows that $K(n)$ has uncountably many $MU$-algebra structures. We first consider the connective Morava $K$-theory spectrum $k(n)$ with $k(n)_*=\bF_p[v_n]$. We choose $x_i$ such that $MU_*=\bZ[x_1,x_2,\ldots]$ and
\[ k(n)=MU/(p,x_1,\ldots,x_{p^n-2},x_{p^n},\ldots).\]
We can also choose these generators in such a way that $x_{p^i-1}$ maps to $v_i$ for $0 \leq i \leq n$ and $x_j$ maps to $0$ otherwise, under a suitable map $MU \to \hE{n}$ (which can be chosen to be $H_\infty$, though it is an open question whether or not it can be chosen to be $E_\infty$).

In this case, $E_1^{0,1}$ is nontrivial, but not large enough to kill all the obstructions. To be more precise, the $E_2$ term for topological Hochschild cohomology of $k(n)$ with ground ring $MU$ looks like
\[ \tilde{E}_2^{*,*} = k(n)_*[\tq_0,\tq_1,\ldots,\tq_{p^n-2},\tq_{p^n},\ldots] \]
with $\tq_i$ in filtration $1$ and total homological degree $-2i-2$.

The term $E_1^{0,1}$ consists of infinite sums $1+\sum v_I Q_I$, where $v_I \in k(n)_*$ is in the appropriate degree and $Q_I=Q_{i_1} \ldots Q_{i_k}$ is a product of Bocksteins. Here $Q_i$ is the Bockstein corresponding to $x_i$ in $MU_*$, or to $\tq_i$ in $\tilde{E}_2^{*,*}$.

Similarly, the term $E_1^{1,1}$ consists of infinite sums $1 \sma 1 + \sum v_{IJ} Q_I \sma Q_J$. The $d_1$ differential $d_1 : E_1^{0,1} \to E_1^{1,1}$ is given by
\[ d_1(v_{ij} Q_i Q_j) = v_{ij} Q_i \sma Q_j - v_{ij} Q_j \sma Q_i,\]
and more generally $d_1(v_I Q_I)$ is given by the sum of all ways to write $I=J \cup K$ of $\pm v_I Q_J \sma Q_K$. In particular, $d_1$ is injective, so $E_2^{0,1}=0$ is trivial.

We also find that $v_{ij} Q_i \sma Q_j = v_{ij} Q_j \sma Q_i$ in $E_2^{1,1}$, and as in \cite[Theorem 3.9]{AnTHH}, the kernel of $d_1 : E_1^{1,1} \to E_1^{2,1}$ picks out the homotopy associative multiplications, and this identifies $E_2^{1,1}$ with $\tilde{E}_2^{2,0}$. Again there can be no $\lim^1$ terms, so the spectral sequence converges completely. This gives a complete description of all the $A_\infty$ structures on $k(n)$ as an $MU$-module. We get the same result for $K(n)$:

\begin{lemma} \label{lem:BRkntoBRKnwe}
The canonical map $\BR^{MU}(k(n)) \to \BR^{MU}(K(n))$ is a weak equivalence.
\end{lemma}

\begin{proof}
This is clear because
\[ \tilde{E}_2^{*,*}(K(n)) \cong v_n^{-1} \tilde{E}_2^{*,*}(k(n)),\]
and these groups are isomorphic in the degrees contributing to $E_2^{*,*}(k(n))$ and $E_2^{*,*}(K(n))$, and the same holds for $E_2^{0,*}$.
\end{proof}

This proves Theorem \ref{thm:main2} for $R=MU$. If $BP$ is a commutative $S$-algebra then the same argument shows that $K(n)$ has uncountably many $BP$-algebra structures.

\subsection{Ground ring $R=S$} \label{ss:BKSSS}
By \cite{AnTHH}, we have an equivalence $THH_{\hE{n}}(K(n)) \to THH_S(K(n))$ (which is visible on $\tilde{E}_2$), and this shows that the $E_2$ term of the spectral sequence converging to $\pi_* \BR^S(K(n))$ is isomorphic to the $E_2$ term of the spectral sequence converging to $\pi_* \BR^{\hE{n}}(K(n))$ in filtration $s \geq 2$. If $p$ is odd this also gives an isomorphism in filtration $s=1$; if $p=2$ there is a possible differential $d_1 : E_1^{0,1} \to E_1^{1,1}$ killing the class $v_n q_{n-1}^2$.

As in \cite{Ra04}, let
\[
\Sigma(n) = K(n)_* \otimes_{BP_*} BP_*BP \otimes_{BP_*} K(n)_*
\cong K(n)_*[t_1,t_2,\ldots]/(v_n t_i^{p^n}-v_n^{p^i} t_i)
\]
be the $n$'th Morava Stabilizer algebra. Here $|t_i|=2(p^i-1)$. Recall \cite{Na02} that, for any choice of multiplication on $K(n)$, we have
\begin{equation*} \label{eq:K(n)_*K(n)}
K(n)_* K(n) \cong \Sigma(n) \otimes \Lambda(\alpha_0,\ldots,\alpha_{n-1})
\end{equation*}
as a ring for $p$ odd, while $\alpha_i^2=t_{i+1}$ for $0 \leq i \leq n-2$ and $\alpha_{n-1}^2=t_n+v_n$ for $p=2$.\footnote{If the reader prefers a unified description of $K(n)_*K(n)$ at all primes it is the above ring with $p$-fold Massey products ($2$-fold Massey products being products) $\langle \alpha_i,\ldots,\alpha_i \rangle = t_{i+1}$ for $0 \leq i \leq n-2$ and $\langle \alpha_{n-1},\ldots,\alpha_{n-1} \rangle = t_n+v_n$ with no indeterminacy.} Here $|\alpha_i|=2p^i-1$.

Also recall that if we consider $K(n)_* K(n)^{op}$ instead we get the same result except that we get to replace the relation $\alpha_{n-1}^2=t_n+v_n$ by $\alpha_{n-1}^2=t_n$ at $p=2$.\footnote{Or replace $\langle \alpha_{n-1},\ldots,\alpha_{n-1} \rangle = t_n+v_n$ with $\langle \alpha_{n-1}, \ldots,\alpha_{n-1} \rangle=v_n$ at any prime.}

We have that
\begin{multline*}
K(n)^*K(n) \cong Hom_{K(n)_*}(K(n)_*K(n),K(n)_*) \\
\cong Hom_{K(n)_*}(\Sigma(n),K(n)_*) \otimes \Lambda(Q_0,\ldots,Q_{n-1}),
\end{multline*}
where $Q_i$ is the Bockstein dual to $\alpha_i$. In particular this means that
\[ E_1^{0,1} \cong \big[ Hom_{K(n)_*}(\Sigma(n),K(n)_*) \otimes \Lambda(Q_0,\ldots,Q_{n-1}) \big]_1^\times, \]
which is large enough to potentially kill all the uniqueness obstructions.

Again there can be no $\lim^1$ terms, so the spectral sequence converges completely. This is clear in positive filtration, for the groups in filtration $0$ this relies on observing that $E_1^{0,t}$ is $p$-torsion.

At $p=2$, a result by Nassau \cite{Na02} gives us our first differential. He shows that if $\phi$ is one multiplication ($A_2$ structure) on $K(n)$ and $\phi^{op}$ is the other, the automorphism $\Xi$ of $K(n)$ given by $t_n \mapsto v_n$ is an antiautomorphism of the multiplication. Hence $\phi$ and $\phi^{op}$ are in the same path component in $\BR^S_2(K(n))$. The difference $\phi-\phi^{op}$ is represented by $v_n q_{n-1}^2$, so $d_1(\Xi)=v_n q_{n-1}^2$.

\section{$S$-algebra $k$-invariants} \label{s:k_invariants}
For connective spectra we can build the $S$-algebra structure by induction on the Postnikov sections. Given a connective spectrum $A$, let $P_m A$ denote the Postnikov section of $A$ with homotopy groups only up to (and including) degree $m$. If $R$ is a connective commutative $S$-algebra then Postnikov sections can be defined in the category of $R$-algebras, so if $A$ is an $R$-algebra then this gives an $R$-algebra structure on $P_m A$ as well. Conversely, $\BR^R(A)=\underset{\longleftarrow}{\lim} \, \BR^R(P_m A)$, so we can understand $\BR^R(A)$ by understanding $\BR^R(P_m A)$ for all $m$.

A theory of $k$-invariants for connective $R$-algebras has been developed by Dugger and Shipley \cite{DuSh}. Suppose $C$ is an $R$-algebra with homotopy groups only up to degree $m-1$, and suppose $M$ is a $\pi_0 C$ module. Let $\cM(C,(M,m))$ be the category of Postnikov extensions of $C$ of type $(M,m)$. The objects are $R$-algebras $Y$ together with a map $Y \to C$ satisfying $\pi_i Y=0$ for $i>m$, $\pi_m Y \cong M$ and $P_{m-1} Y \simeq C$. The morphism are maps over $C$ inducing an isomorphism on homotopy.

\begin{theorem} (Dugger-Shipley, \cite[Proposition 1.5]{DuSh}) With $\mathcal{M}(C,(M,n))$ as above,
\begin{equation*}
\pi_0 \cM(C,(M,m)) \cong THH^{m+2}_R(C;M)/Aut(M).
\end{equation*}
\end{theorem}

Now suppose $C=P_{m-1} A$, $M=\pi_m A$, and we want to make sure that $Y \in \cM(C,(M,m))$ has the homotopy type of $P_m A$. Then $Y$ has to have the correct additive $k$-invariant, which is a map $C \to \Sigma^{m+1} HM$. Recall that the topological Hochschild cohomology spectral sequence converging to $THH^*_R(C;M)$ has $\tilde{E}_1^{s,t}=[\Sigma^t C^{(s)},M]$, contributing to $\pi_{t-s} THH_R(C;M)=THH^{s-t}_R(C;M)$. In particular, the additive $k$-invariant of $Y$ is an element in $\tilde{E}_1^{1,-m-1}$, contributing to $THH^{m+2}_R(C;M)$.

If the additive $k$-invariant $k_m$ is trivial then $Y \simeq C \vee \Sigma^m HM$ as a spectrum, and $Y$ always has at least one $S$-algebra structure, namely the square zero extension. If $k_m$ is non-trivial, it might or might not survive the topological Hochschild cohomology spectral sequence. If $d_r(k_m)=y \neq 0$ then $y$ represents the obstruction to extending the $S$-algebra structure on $C$ to an $S$-algebra structure on $Y$. If $k_m$ survives then $Y$ has at least one $S$-algebra structure.

\section{$k$-invariants for Morava $K$-theory} \label{s:kforkn}
Again we study the moduli problem over each ground ring separately. First we use $BP \langle n \rangle_p$, which has homotopy groups
\[ (BP \langle n \rangle_p)_* = \bZ_p[v_1,\ldots,v_n] \]
and is the appropriate connetive version of $\hE{n}$, as the ground ring, assuming it can be given the structure of a commutative $S$-algebra. Then we use $MU$, and finally we use the sphere spectrum $S$. First we recall the following change-of-rings result:

\begin{lemma}[{\cite[Corollary 2.5]{AHL}}] \label{lemma:changeofcoeff}
Suppose $A \to B$ is a map of $S$-algebras and $M$ is an $A-B$-bimodule, given an $A-A$-bimodule structure by pullback. Then there is a spectral sequence
\[ \tilde{E}_2^{*,*}=Ext_{\pi_* B \sma_R A}^{**}(B_*,M_*) \Longrightarrow \pi_* THH_R(A;M).\]
\end{lemma}
In particular, when $B=M=H\bF_p$ we get a spectral sequence
\[ \tilde{E}_2^{*,*}=Ext_{H_*^R(A;\bF_p)}(\bF_p,\bF_p) \Longrightarrow \pi_* THH_R(A;H\bF_p),\]
where $H_*^R(A;\bF_p)$ denotes $\pi_* A \sma_R H\bF_p$.

Since $k(n)$ has homotopy in degrees that are multiplies of $2p^n-2$, let $q=2p^n-2$. Each additive $k$-invariant $k_m \in H_R^{mq+1}(P_{(m-1)q} k(n);\bF_p)$ is nontrivial, this follows by considering $H_*(P_{mq} k(n);\bF_p)$, which is different from $H_*(P_{(m-1)q} k(n) \vee \Sigma^{mq} H\bF_p;\bF_p)$.

\subsection{Ground ring $BP \langle n \rangle_p$} \label{ss:k_invBPn}
Since we are planning to use Lemma \ref{lemma:changeofcoeff}, we start by calculating the $\BP{n}_p$ homology of the Postnikov sections of $k(n)$.

\begin{prop}
The $\BP{n}_p$-homology of $H\bF_p$, $P_{mq} k(n)$ and $k(n)$ is as follows:
\begin{enumerate}
\item $H_*^{\BP{n}_p}(H\bF_p;\bF_p)=\Lambda_{\bF_p}(\alpha_0,\ldots,\alpha_n)$,
\item $H_*^{\BP{n}_p}(P_{mq} k(n);\bF_p)=\Lambda_{\bF_p}(\alpha_0,\ldots,\alpha_{n-1},a_{m+1})$,
\item $H_*^{\BP{n}_p}(k(n);\bF_p)=\Lambda_{\bF_p}(\alpha_0,\ldots,\alpha_{n-1})$.
\end{enumerate}
Here $\alpha_i$ is in degree $2p^i-1$ and $a_{m+1}$ is in degree $(m+1)q+1$, $a_1=\alpha_n$.
\end{prop}

\begin{proof}
This is clear, using that we can write
\begin{eqnarray*}
H\bF_p & = & \BP{n}_p/(p,v_1,\ldots,v_n), \\
P_{mq} k(n) & = & \BP{n}_p/(p,v_1,\ldots,v_{n-1},v_n^{m+1}), \\
k(n) & = & \BP{n}/(p,v_1,\ldots,v_{n-1}).
\end{eqnarray*}
\end{proof}

\begin{prop}
Assuming that $\BP{n}_p$ is a commutative $S$-algebra, topological Hochschild cohomology of $H\bF_p$, $P_{mq} k(n)$ and $k(n)$ over $\BP{n}_p$ with coefficients in $H \bF_p$ is as follows:
\begin{enumerate}
\item $THH^*_{\BP{n}_p}(H\bF_p; H\bF_p) \cong \bF_p[q_0,\ldots,q_n]$,
\item $THH^*_{\BP{n}_p}(P_{mq} k(n); H\bF_p) \cong \bF_p[q_0,\ldots,q_{n-1}, b_{m+1}]$,
\item $THH^*_{\BP{n}_p}(k(n); H\bF_p) \cong \bF_p[q_0,\ldots,q_{n-1}]$.
\end{enumerate}
Here $q_i$ is in cohomological degree $2p^i$ and $b_{m+1}$ is in degree $(m+1)q+2$, $b_1=q_n$.
\end{prop}

\begin{proof}
We use Lemma \ref{lemma:changeofcoeff}. In each case there can be no differentials, because the $E_2$ term is concentrated in even total degree.
\end{proof}

The additive $k$-invariant of $k(n)$ dictates that we choose the $k$-invariant in
\[ THH^{mq+2}_{\BP{n}_p}(P_{(m-1)q} k(n);H\bF_p) \]
as $b_m+f(q_0,\ldots,q_{n-1})$ where $f$ has degree at least $2$ in the $q_i$'s.

Next we compare this with the moduli space of $\hE{n}$-algebra structures on $K(n)$.

\begin{lemma} \label{lem:fowdk}
Assuming that $\BP{n}_p$ is a commutative $S$-algebra, the canonical maps
\[ \BR^{\BP{n}_p}(k(n)) \to \BR^{\BP{n}_p}(K(n)) \to \BR^{\hE{n}}(K(n)) \]
are weak equivalences.
\end{lemma}

\begin{proof}
This is similar to the proof of Lemma \ref{lem:BRkntoBRKnwe}
\end{proof}

Now we can compare the two methods of studying the set of equivalence classes of $\BP{n}_p$-algebra structures on $k(n)$. We find that in the spectral sequence converging to $\pi_* \BR^{\BP{n}_p}(k(n))$, each uniqueness obstruction is represented by a class
\[ v_n^m f(q_0,\ldots,q_{n-1}) \]
for some $m \geq 1$, where $f(q_0,\ldots,q_{n-1})$ has homological degree $-mq-2$. If $f(q_0,\ldots,q_{n-1})=q_{i_1} \cdots q_{i_j}$ has degree $j$ in the $q_i$'s this represents changing the $A_j$ structure by the map $\Sigma^{j-2} k(n)^{(j)} \to k(n)$ given by first applying
\[ Q_f = Q_{i_1} \sma \ldots \sma Q_{i_j} \]
and then multiplying the factors and multiplying by $v_n^m$.

On the other hand, we can interpret the polynomial $f(q_0,\ldots,q_{n-1})$ as being an element of $THH^{mq+2}(P_{(m-1)q} k(n); H\bF_p)$, represented in the topological Hochschild cohomology spectral sequence by the composite
\[ (P_{(m-1)q} k(n))^{(j)} \overset{Q_f}{\to} \Sigma^{mq-j+2} (P_{(m-1)q} k(n))^{(j)} \to \Sigma^{mq-j+2} H\bF_p.\]

\begin{lemma}
Given a uniqueness obstruction $v_n^m f(q_0,\ldots,q_{n-1})$ of degree $j$ in the $q_i$'s represented by $v_n^m Q_f : \Sigma^{j-2} k(n)^{(j)} \to k(n)$, we get a commutative diagram as follows:
\[ \xymatrix{ & \Sigma^{j-2} (P_{(m-1)q} k(n))^{(j)} \ar[r]^-{Q_f} & \Sigma^{mq} H\bF_p \ar[d] \\
\Sigma^{j-2} k(n)^{(j)} \ar[r]^-{v_n^m Q_f} \ar[ur] & k(n) \ar[r] & P_{mq} k(n) \ar[d] \\ & & P_{(m-1)q} k(n) } \]
\end{lemma}

\begin{proof}
Consider the following commutative diagram:
\[ \xymatrix{ \Sigma^{j-2} (P_{mq} k(n))^{(j)} \ar[r]^-{Q_f} \ar[d] & \Sigma^{mq} P_{mq} k(n) \ar[d] \ar[r]^-{v_n^m} & P_{mq} k(n) \ar[d] \\
\Sigma^{j-2} (P_{(m-1)q} k(n))^{(j)} \ar[r]^-{Q_f} & \Sigma^{mq} P_{(m-1)q} k(n) \ar[r]^-{v_n^m=0} \ar@{.>}[ur] & P_{(m-1)q} k(n) } \]
This gives us a map $\Sigma^{j-2} (P_{(m-1)q} k(n))^{(j)} \to P_{mq} k(n)$, and this map is trivial on $P_{(m-1)q} k(n)$ so it factors through $\Sigma^{mq} H\bF_p$.
\end{proof}

The upshot of this is that we can translate from obstructions in the spectral sequence converging to $\pi_* \BR^{\BP{n}}(k(n))$, which by Lemma \ref{lem:fowdk} and the equivalence between topological Hochschild cohomology over $\hE{n}$ and $S$ are the obstructions in the spectral sequence converging to $\pi_* \BR^S(K(n))$ to obstructions in the spectral sequence converging to $THH_{\BP{n}}^*(P_{(m-1)q} k(n);H\bF_p)$.

\subsection{Ground ring $MU$} \label{ss:k_invMU}
Next we do the same with $MU$ as the ground ring. If we knew that $\BP{n}_p$ was a commutative $S$-algebra then this section would not be necessary. The corresponding results are as follows:

\begin{prop}
The $MU$-homology of $H\bF_p$, $P_{mq} k(n)$ and $k(n)$ is as follows:
\begin{enumerate}
\item $H_*^{MU}(H\bF_p;\bF_p) \cong \Lambda_{\bF_p}(\aalpha_0,\aalpha_1,\ldots)$,
\item $H_*^{MU}(P_{mq} k(n);\bF_p) \cong \Lambda_{\bF_p}(\aalpha_i \, : \, i \neq p^n-1, a_{m+1})$,
\item $H_*^{MU}(k(n);\bF_p) \cong \Lambda_{\bF_p}(\aalpha_i \, : \, i \neq p^n-1)$.
\end{enumerate}
Here $\aalpha_i$ is in degree $2i+1$ and $a_{m+1}$ is in degree $(m+1)q+1$, $a_1=\aalpha_{p^n-1}$.
\end{prop}

\begin{proof}
This is clear, using that we can write
\begin{eqnarray*}
H\bF_p & = & MU/(p,x_1,\ldots,x_{p^n-2}, x_{p^n-1}, x_{p^n}, \ldots), \\
P_{mq} k(n) & = & MU/(p,x_1,\ldots,x_{p^n-2}, x_{p^n-1}^m, x_{p^n},\ldots), \\
k(n) & = & MU/(p,x_1,\ldots,x_{p^n-2}, x_{p^n}, \ldots).
\end{eqnarray*}
\end{proof}

Now we can calculate topological Hochschild cohomology:

\begin{prop}
Topological Hochschild cohomology of $\bF_p$, $P_{mq} k(n)$ and $k(n)$ over $MU$ with coefficients in $H\bF_p$ is as follows:
\begin{enumerate}
\item $THH_{MU}^*(H\bF_p;H\bF_p) \cong \bF_p[\tq_0,\tq_1,\ldots]$,
\item $THH_{MU}^*(P_{mq} k(n);H\bF_p) \cong \bF_p[\tq_i \, : \, i \neq p^n-1, b_{m+1}]$,
\item $THH_{MU}^*(k(n);H\bF_p) \cong \bF_p[\tq_i \, : \, i \neq p^n-1]$.
\end{enumerate}
Here $\tq_i$ is in cohomological degree $2i+2$ and $b_{m+1}$ is in degree $(m+1)q+2$, $b_1=\tq_{p^n-1}$.
\end{prop}

\begin{proof}
Again this follows from Lemma \ref{lemma:changeofcoeff}.
\end{proof}

Recall from Lemma \ref{lem:BRkntoBRKnwe} that $\BR^{MU}(k(n)) \to \BR^{MU}(K(n))$ is a weak equivalence. Just as with $\BP{n}_p$ as the ground ring, we can translate from obstructions in the spectral sequence converging to $\BR^S(K(n))$ to obstructions in the spectral sequence converging to $THH^*_{MU}(P_{(m-1)q} k(n);H\bF_p)$. In this case, only $\bF_p[\tq_0, \tq_{p-1},\ldots,\tq_{p^{n-1}-1}]$ correspond to obstructions in the spectral sequence converging to $\pi_* \BR^S(K(n))$.

By this we mean that the $MU$-algebra $k$-invariant for building $P_{mq} k(n)$ from $P_{(m-1)q} k(n)$ lives in $THH^{mq+2}(P_{(m-1)q} k(n);H\bF_p)$ and looks like $b_m + f(\tq_i \, : \, i \neq p^n-1)$ where $f$ has degree at least $2$ in the $\tq_i$'s. This corresponds to the uniqueness obstruction $v_n^m f(\tq_i \, : \, i \neq p^n-1)$ in the $E_2$-term of the spectral sequence converging to $\pi_* \BR^{MU}(K(n))$, and the canonical map $\BR^{MU}(K(n)) \to \BR^S(K(n))$ induces a map on $E_2$-terms, under which $\tq_{p^i-1}$ maps to $q_i$. This is clear, because both $\tq_{p^i-1}$ and $q_i$ are represented by the Bockstein corresponding to $v_i$.

\subsection{Ground ring $S$} \label{ss:k_invS}
Finally we do the same with $S$ as the ground ring. Let $\bar{A}_*$ denote the dual Steenrod algebra with $\ttau_n$ missing, or with $\xxi_{n+1}$ missing but with $\xxi_{n+1}^2$ present if $p=2$. In the following we will state all results at odd primes and leave the standard modifications, replacing $\ttau_i$ with $\xxi_i$ and $\xxi_i$ with $\xxi_i^2$ at $p=2$ to the reader.

\begin{prop}
The mod $p$ homology of $\bF_p$, $P_{mq} k(n)$ and $k(n)$ is as follows:
\begin{enumerate}
\item $H_*(\bF_p;\bF_p) \cong A_*$,
\item $H_*(P_{mq} k(n);\bF_p) \cong \bar{A}_* \otimes \Lambda_{\bF_p}(a_{m+1})$,
\item $H_*(k(n);\bF_p) \cong \bar{A}_*$.
\end{enumerate}
Here $a_{m+1}$ is in degree $(m+1)q+1$, $a_1=\ttau_n$.
\end{prop}

\begin{proof}
Only part $2$ is not well known. Consider the long exact sequence obtained by taking the mod $p$ homology of the (co)fiber sequence
\[ \Sigma^{mq} H\bF_p \to P_{mq} k(n) \to P_{(m-1)q} k(n) \to \Sigma^{mq+1} H\bF_p. \]
By induction we have $H_*(P_{(m-1)q} k(n);\bF_p) \cong \bar{A}_* \otimes \Lambda_{\bF_p}(a_m)$, and the map to $H_*(\Sigma^{mq+1} H\bF_p; \bF_p)$ is determined by being $\bar{A}_*$-linear and that $1 \mapsto 0$ and $a_m \mapsto \Sigma^{mq+1} 1$. The result follows by combining the kernel and cokernel of this map.
\end{proof}

\begin{thm}
Topological Hochchild cohomology of $H\bF_p$, $P_{mq} k(n)$ and $k(n)$ with coefficients in $H\bF_p$ is as follows:
\begin{enumerate}
\item $THH_S^*(H\bF_p;H\bF_p) \cong P_p(\delta \ttau_0,\delta \ttau_1,\ldots)$,
\item $THH_S^*(P_{mq} k(n);H\bF_p) \cong \Lambda(\delta \xxi_{n+1}) \otimes P_p(\delta \ttau_i \, : \, i \neq n) \otimes \bF_p[b_{m+1}]$,
\item $THH_S^*(k(n);H\bF_p) \cong \Lambda(\delta \xxi_{n+1}) \otimes P_p(\delta \ttau_i \, : \, i \neq n)$.
\end{enumerate}
\end{thm}

\begin{proof}
The first part is dual to B\"okstedt's original calculation of topological Hoch\-schild homology of $\bF_p$ \cite{Bo2}. For 2, consider the spectral sequence
\[ E_2 = \Lambda(\delta \xxi_i \, : \, i \geq 1) \otimes \bF_p[\delta \ttau_i \, : \, i \neq n] \otimes \bF_p[b_{m+1}] \Longrightarrow THH_S^*(P_{mq} k(n);\bF_p) \]
from Lemma \ref{lemma:changeofcoeff}. The map $P_{mq} k(n) \to H\bF_p$ induces a map on topological Hoch\-schild cohomology in the opposite direction, inducing differentials $d_{p-1}(\delta \xxi_{i+1}) = (\delta \ttau_i)^p$ for $i \neq n$.

The class $b_{m+1}$ is the next additive $k$-invariant for $k(n)$, and because we know that $k(n)$ can be given an $S$-algebra structure, $b_{m+1}$ has to survive the spectral sequence. The class $\delta \xxi_{n+1}$ survives for degree reasons, so each generator is a permanent cycle. Using the multiplicative structure, the spectral sequence collapses at the $E_p$ term and part 2 of the theorem follows. Part 3 is similar.
\end{proof}

We note that the $p$'th powers of $\delta \ttau_i$ for $0 \leq i \leq n-1$ all die, and we make the following simple but cruical observation:

\begin{lemma} \label{lem:noobinrightdeg}
Consider the $k$-invariant for $k(n)$ in $THH_S^{mq+2}(P_{(m-1)q} k(n); \bF_p)$. There are no polynomials $f(\delta \ttau_0,\ldots,\delta \ttau_{n-1}) \in P_p(\delta \ttau_0,\ldots,\delta \ttau_{n-1})$ in this degree.
\end{lemma}

\begin{proof}
This is clear because the element in highest degree is $(\delta \ttau_0)^{p-1} \cdots (\delta \ttau_{n-1})^{p-1}$ in degree $2p^n-2$, which is less than $mq+2$.
\end{proof}

Of course the generators $\delta \ttau_i$ are related to the generators $q_i$ and $\tq_j$ from the previous sections:

\begin{lemma} \label{lem:gensarethesame}
The canonical map
\[ THH_{MU}^*(P_{mq} k(n);\bF_p) \to THH_S^*(P_{mq} k(n);\bF_p), \]
maps $\tq_j$ to $\delta \ttau_i$ if $p^i-1=j$ and $0$ otherwise.

Similarly, if $\BP{n}_p$ is a commutative $S$-algebra, the canonical map
\[ THH_{\BP{n}_p}^*(P_{mq} k(n);\bF_p) \to THH_S^*(P_{mq} k(n);\bF_p) \]
maps $q_i$ to $\delta \ttau_i$.
\end{lemma}

\begin{proof}
This follows by the description of all of the $E_2$-terms in terms of Bocksteins.
\end{proof}

\section{Proof of Theorem \ref{thm:main1}} \label{s:proof}
We are now in a position to prove Theorem \ref{thm:main1}. As we have seen, each uniqueness obstruction looks like $v_n^m f(q_0,\ldots,q_{n-1})$ for some $m \geq 1$ and some monomial $f(q_0,\ldots,q_{n-1})$, and we can find these uniqueness obstructions in the corresponding spectral sequence converging to $THH_{MU}^*(P_{(m-1)q} k(n);H\bF_p)$.

In the corresponding spectral sequence converging to $THH_S^*(P_{(m-1)q} k(n);\bF_p)$, $f(q_0,\ldots,q_{n-1})$ is killed by a differential, which means that the corresponding $S$-algebra structures on $P_{mq} k(n)$ are equivalent. By considering the pullback square
\[ \xymatrix{ PB \simeq k(n) \ar[d] \ar[r] & k(n) \ar[d] \\ P_{mq} k(n) \ar[r] & P_{mq} k(n) } \]
of $S$-algebras, we see that the equivalence can be lifted to $k(n)$. Now we can invert $v_n$ by $K(n)$-localizing, so this gives an equivalence between the corresponding $S$-algebra structures on $K(n)$ as well.

We claim that this is enough to conclude that the obstructions are also killed in the spectral sequence converging to $\pi_* \BR^S(K(n))$.  To see this, consider $k(n)$ and $K(n)$ as $MU$-modules, and consider the following commutative diagram:
\[ \xymatrix{ \BR^{MU}(k(n)) \ar[r] \ar[d]^\simeq & \BR^S(k(n)) \ar[d] \\
\BR^{MU}(K(n)) \ar[r] & \BR^S(K(n)) } \]
We showed in Lemma \ref{lem:BRkntoBRKnwe} that $\BR^{MU}(k(n)) \to \BR^{MU}(K(n))$ is a weak equivalence, and we understand the $E_2$ terms of the spectral sequences converging to the homotopy groups of all the spaces in the diagram except for $\BR^S(k(n))$. The spectral sequences converging to $\pi_* \BR^{MU}(k(n))$ and $\pi_* \BR^{MU}(K(n))$ collapse, and from the $E_2$ terms we can read off that the map $\pi_0 \BR^{MU}(k(n)) \to \pi_0 \BR^S(K(n))$ is surjective. 

In $\pi_0 \BR^{MU}(k(n))$, there are classes that map surjectively onto the $E_2$ term of the spectral sequence converging to $\pi_* \BR^S(K(n))$ which are all hit by differentials in the spectral sequence converging to $THH_S^*(P_{(m-1)q} k(n); H\bF_p)$ for some $m$ (Lemma \ref{lem:noobinrightdeg} and \ref{lem:gensarethesame}), hence the same must happen in the spectral sequence converging to $\pi_* \BR^S(K(n))$.

Our argument would be simplified by the existence of a commutative $S$-algebra structure on $\BP{n}_p$, in which case it follows that all the uniqueness obstructions for building $k(n)$ as a $BP \langle n \rangle_p$-algebra are hit by differentials in the spectral sequence converging to $THH_S^*(P_{(m-1)q} k(n);H\bF_p)$ for some $m$. In particular, when $n=1$ using $\ell_p$ instead of $MU$ gives a simpler argument.

\section{$2$-periodic Morava $K$-theory} \label{s:2periodic}
There is a $2$-periodic version of Morava $K$-theory, given by
\[ K_n=E_n/(p,u_1,\ldots,u_{n-1}), \]
where $E_n$ is the Morava $E$-theory spectrum associated to a formal group of height $n$ over a perfect field $k$ of characteristic $p$. The spectrum $E_n$ is a commutative $S$-algebra \cite{GoHo}, and $K_n$ has homotopy groups
\[ (K_n)_* \cong k[u,u^{-1}] \]
with $|u|=2$. We can also ask about the space of $S$-algebra structures on $K_n$. When $p=2$ and $n=1$, $K_n=K(n)$, if $p>2$ or $n>1$ the author \cite{AnTHH} found that $THH_S(K_n)$ varies over the moduli space of $S$-algebra structures on $K_n$, so there can be no unique $S$-algebra structure on $K_n$.

\begin{conj}
There are only finitely many $S$-algebra structures on $K_n$, in the sense that the moduli space of $S$-algebra structures on $K_n$ has finitely many components.
\end{conj}

\begin{proof}[Outline of possible proof]
The spectral sequence converging to $\pi_* \BR^S(K_n)$ is very similar to the one converging to $\pi_* \BR^S(K(n))$, but now each of the $n$ polynomial generators are in degree $-2$ instead of degree $-2p^i$ for $0 \leq i \leq n-1$.

If we try to build the connective version $k_n$ using its Postnikov tower we need to understand the topological Hochschild cohomology spectral sequence. Since $H_*(k_n;\bF_p) \cong H_*(k(n);\bF_p) \otimes P_{p^n-1}(u)$, and similarly for the Postnikov sections, we get some extra classes in the $E_2$ term. Assuming that these classes are permanent cycles, we find that we have more choices than before. To build $P_2 k_n$ from $Hk$, we need a class in
\[ THH_S^4(Hk;Hk),\]
and for $p$ odd we are free to choose $(\delta \ttau_0)^2$. If $p=2$ and $n>1$ we can choose $\delta \xxi_2$. In each case this corresponds to a noncommutative multiplication. Next, to build $P_4 k_n$ from $P_2 k_n$ we need a class in $THH^6_S(P_2 k_n;Hk)$. If $p>3$ we can choose the class we need for $u$ to square to something nontrivial plus $(\delta \ttau_0)^3$, if $n \geq 2$ and $p=2$ or $p=3$ there are similar choices.

However, assuming that the additional classes do not change the behavior of the spectral sequence, the $p$'th powers of $\delta \ttau_0, \ldots, \delta \ttau_{n-1}$ still die, so for $m$ sufficiently large there are no such classes in $THH_S^{2m+2}(P_{2m-2} k_n;Hk)$.
\end{proof}

\bibliographystyle{plain}
\bibliography{b}

\end{document}